\documentclass[12pt,a4paper]{article}
\usepackage[latin2]{inputenc}
\usepackage{amsmath}
\usepackage{amsfonts}
\usepackage{amssymb}
\usepackage{graphicx}
\usepackage[left=2.00cm, right=2.00cm, top=2.50cm, bottom=2.50cm]{geometry}
\usepackage[T1]{fontenc}

\def\Q{{\mathbb Q}}
\def\Z{{\mathbb Z}}

\def\ov{\overline}

\newtheorem{lemma}{Lemma}

\title{
Calculating generators of power integral bases \\
in sextic fields with a real quadratic subfield
}
\author{
Istv\'an Ga\'al\\
{\small University of Debrecen, Mathematical Institute} \\
{\small H--4002 Debrecen Pf.400., Hungary,} \\
{\small e--mail: gaal.istvan@unideb.hu},
}

\begin{document}
\baselineskip=17pt

\maketitle
\thispagestyle{empty}

\renewcommand{\thefootnote}{\arabic{footnote}}
\setcounter{footnote}{0}

\vspace{0.5cm}


\begin{abstract}
We discuss the problem of calculating generators of power integral bases
in sextic fields, especially 
focusing on the case of sextic fields with real quadratic subfields.
Our main purpose is to describe an efficient algorithm
for calculating generators of power integral bases. 
We show that appropriately using integer arithmetics 
speeds up the calculations  considerably.
Our experiences lead to some interesting general statements
on generators of power integral bases in number fields
generated by a unit.
\end{abstract}

\section{Introduction}

Monogenity and power integral bases is an important classical 
area of algebraic number theory going back to 
Dedekind \cite{dedekind}, Hasse  \cite{Ha} and Hensel \cite{He}, cf. also \cite{nark}.
For the present state of this area we refer to \cite{book}. 
It is called a problem of Hasse to give an arithmetic characterization of those number fields which have a power integral basis.

We discuss the problem of calculating
generators of power integral bases in sextic fields with a quadratic subfield,
focusing for the case of real quadratic subfields.
Our main result is to give an efficient algorithm 
for calculating generators of power integral bases
in sextic fields with a real quadratic subfield.

In Section \ref{s1} we present the related notions, 
discuss some former results 
and formulate some auxiliary statements we apply.

In Section \ref{s2} we discuss an algorithm which turns out to be 
faster in our case, than the previously used methods.

Finally, in Section \ref{calc} we provide some examples for the
application of our algorithm. The numerical results of our calculations
lead us to recognize some general properties of generators of power integral
bases: these simple but useful properties are also given in this section.

\section{Preliminaries}
\label{s1}

A number field $K$ of degree $n$ with a ring of integers $\Z_K$ 
is called {\it monogenic} (cf. \cite{book}) if there 
exists $\xi\in \Z_K$ such that $(1,\xi,\ldots,\xi^{n-1})$ is an integral basis, 
called {\it power integral basis}. 
We call $\xi$ the {\it generator} of this power integral basis.
$\alpha,\beta\in\Z_K$ are called {\it equivalent},
if $\alpha+\beta\in\Z$ or $\alpha-\beta\in\Z_K$. 
Obviously, $\alpha$ generates a power integral basis in $K$
if and only if any $\beta$, equivalent to $\alpha$, does.
As it is known, any algebraic number field admits up to equivalence only
finitely many generators of power integral bases.

An irreducible polynomial $f(x)\in\Z[x]$ is called {\it monogenic}, 
if a root $\xi$ of $f(x)$ generates a power integral basis in $K=\Q(\xi)$.
If $f(x)$ is monogenic, then $K$ is also monogenic, but the converse is
not true.

The field $K$ is called {\it relative monogenic} over its subfield $M$, 
if $K$ has a relative integral basis over $M$ of type 
$(1,\gamma,\ldots,\gamma^k)$ where $k=[K:M]$.

For $\alpha\in\Z_K$ (generating $K$ over $\Q$)  the module index 
\[
I(\alpha)=(\Z_K:\Z[\alpha])
\]
is called the {\it index} of $\alpha$. The element 
$\alpha$ generates a power integral basis in $K$
if and only if $I(\alpha)=1$. If $\alpha^{(i)}\; (1\le i\le n)$ 
are the conjugates of 
$\alpha$ in $K$  of degree $n$, then
\[
I(\alpha)=\frac{1}{\sqrt{|D_K|}}\prod_{1\le i<j\le n}|\alpha^{(i)}-\alpha^{(j)}|,
\]
where $D_K$ is the discriminant of $K$.
For more details concerning 
the classical and also very actual topics of  monogenity 
and power integral bases cf. \cite{book}.

The powerful methods of Newton polygons and Dedekind criterion
are intensively used during the last couple of years.
These can usually prove the non-monogenity and in certain cases the monogenity 
of  number fields, see \cite{axioms}.
On the other hand it is also an important problem to determine all non-equivalent generators of power integral bases. 
For this purpose we need to use Diophantine methods, 
solving the index form equation corresponding to an integral basis
(cf. \cite{book}).
There are efficient algorithms for cubic and quartic number fields,
but only partial results for higher degree number fields, like
sextic fields with a quadratic subfield or octic fields with a 
quadratic subfield.
These algorithms for the "complete resolution" of index form equations,
may require too long CPU time. On the other hand, it turned out, 
that often there are some
fast algorithms that produce all generators of power integral bases, with 
coefficients, say $<10^{100}$ in absolute value, with respect to an integral basis.
All experiences show, that generators of power integral bases 
have very small coefficients in an integral basis, hence these algorithms
give all generators with a high probability, certainly all generators that
can be used in practice or for further calculations. For such algorithms see e.g.
\cite{g24f}.

In this paper we consider {\it sextic fields}. As we have seen, the algorithm \cite{bilu}
for general sextic fields requires a huge amount of CPU time, therefore
not feasible to apply. 

A couple of results discuss the easy case, when
the sextic field is a {\it composite of a quadratic and a cubic subfield}. 
In these cases the relative index over the quadratic subfield 
provides a cubic relative Thue equation and we have two further factors of
the index, which, together allow to determine generators of power integral bases,
cf. e.g. \cite{grs}.

In case of sextic fields with a
{\it complex quadratic subfield} $M$, the fast algorithm \cite{relthue}
for solving the relative Thue equation appearing in the calculation
is very efficient,
and $M$ has only a few trivial units, which makes the procedure
much easier, see \cite{gp1996}, \cite{grs62}.

Therefore we focus in this paper for sextic fields with 
a real quadratic subfield.

\vspace{1cm}

In our algorithm we shall use some consequences
of previous results, which we present here for completeness.
The following Lemmas are special cases of the general results 
of \cite{grs}, formulated for our case of a sextic field
$K$ with a quadratic subfield $M$ with a ring of integers $\Z_M$.

\begin{lemma}\mbox{} \\
A. If $K$ is monogenic, then $K$ is also relative monogenic over $M$.\\
B. All generators of power integral bases of $K$ are of the form
\[
\gamma=A+\nu \gamma_0,
\]
where $A\in\Z_M$, $\nu$ is a unit in $M$ and $\gamma_0$
generates a relative power integral basis of $K$ over $M$.
\label{th2}
\end{lemma}

Assume $K=M(\alpha)$ and $f(x)=x^3+f_2x^2+f_1x+f_0\in \Z_M$
is  the relative defining polynomial of $\alpha$ over $M$. 
Denote by $\beta=\beta^{(1)},\beta^{(2)} $ the conjugates of any $\beta\in M$.
Accordingly, $f^{(i)}(x)=x^3+f_2^{(i)}x^2+f_1^{(i)}x+f_0^{(i)}$ 
are the conjugates of $f(x)$ ($i=1,2$).
Let $\alpha^{(i,j)}$ be the roots of $f^{(i)}(x)$ ($i=1,2,j=1,2,3$) 
and denote the conjugates of any 
$\gamma\in K$ similarly.

\begin{lemma} For  $\gamma\in\Z_K$ generating $K$ over $\Q$ we have
\[
I(\gamma)=I_{K/M}(\gamma)\cdot J(\gamma)
\]
where
\[
I_{K/M}(\gamma)=(\Z_K:\Z_M[\gamma])=
\frac{1}{\sqrt{|N_{M/\Q}(D_{K/M})|}}
\prod_{i=1}^2\prod_{1\le j_1< j_2\le 3}|\gamma^{(i,j_1)}-\gamma^{(i,j_2)}|
\]
is the relative index of $\gamma$ and
\[
J(\gamma)=\frac{1}{|D_M|^{3/2}}
\prod_{j_1=1}^3\prod_{j_2=1}^3 |\gamma^{(1,j_1)}-\gamma^{(2,j_2)}|.
\]
\label{relindex}
\end{lemma}

Note that if $I(\gamma)=1$, that is $\gamma$ generates a power integral basis 
of $K$, then $I_{K/M}(\gamma)=1$ and $J(\gamma)=1$.

\section{The algorithm}
\label{s2}

In this section we formulate exactly our computational task,
discuss the steps of our calculations and emphasize 
how the application of integer arithmetics speeds up the calculations
considerably.

Let $M=\Q(\sqrt{m})$  ($m>1$ square free) be a real quadratic field.
Let $\omega=\sqrt{m}$ if $m\equiv 2,3\; (\bmod \; 4)$ and 
$\omega=(1+\sqrt{m})/2$ if $m\equiv 1\; (\bmod \; 4)$.
Let $K=M(\alpha)$, where $f(x)=x^3+f_2x^2+f_1x+f_0\in \Z_M$
is  the relative defining polynomial of $\alpha$ over $M$, as above. 

For simplicity's sake assume that we can choose the generator element 
$\alpha$ of $K$ so that $(1,\alpha,\alpha^2,\omega\alpha,\omega\alpha^2)$ 
is an integral basis of $K$ (according to the tables
of  \cite{BMO90}, \cite{Oli89}  in about 99\% 
of the cases $\alpha$ can be chosen with this property, see \cite{book}).

Then any $\gamma\in \Z_K$ can be written in the form
\begin{equation}
\gamma=A+X\alpha+Y\alpha^2,
\label{AXY}
\end{equation}
with $A,X,Y\in\Z_M$.

Set
$A=a_1+\omega a_2, X=x_1+\omega x_2,Y=y_1+\omega y_2$ with $a_1,a_2,x_1,x_2,y_1,y_2\in\Z$ in the representation (\ref{AXY}).
To determine $\gamma$ up to equivalence, we need to calculate
$a_2,x_1,x_2,y_1,y_2$.
\\

\noindent
{\bf The purpose of our calculations}\\
Assume our task is to determine all generators $\gamma$ 
of power integral bases of $K$ with 
\begin{equation}
\max(|a_2|,|x_1|,|x_2|,|y_1|,|y_2|)<C,
\label{bound}
\end{equation}
with $C=10^{50}$, say. According to all experiences, these
generators provide all possible generators of power integral bases with
high probability.
\\

\noindent
{\bf Step 1}\\
First we calculate generators of relative power integral bases of $K$ 
over $M$.
We have
\[
\prod_{i=1}^2\prod_{1\le j_1\le j_2\le 3}|\gamma^{(i,j_1)}-\gamma^{(i,j_2)}|
=\prod_{i=1}^2\prod_{1\le j_1\le j_2\le 3}
|\alpha^{(i,j_1)}-\alpha^{(i,j_2)}|\cdot
|X+(\alpha^{(i,j_1)}+\alpha^{(i,j_2)})Y|.
\]
Using the quadratic coefficient $f_2$ of $f(x)$, for $\{j_1,j_2,j_3\}=\{1,2,3\}$ we have
$f_2^{(i)}=-\alpha^{(i,j_1)}-\alpha^{(i,j_2)}-\alpha^{(i,j_3)}$, whence
by
\[
\sqrt{|N_{M/\Q}(D_{K/M})|}=\prod_{i=1}^2\prod_{1\le j_1\le j_2\le 3}
|\alpha^{(i,j_1)}-\alpha^{(i,j_2)}|
\]
we obtain
\[
I_{K/M}(\gamma)=N_{M/\Q}(N_{K/M}(X-\vartheta Y)),
\]
where $\vartheta=f_2+\alpha$.

Hence in view of Lemma \ref{relindex}
the relative index of $\gamma$ gives rise to the cubic relative Thue equation
\begin{equation}
N_{K/M}(X-\vartheta Y)=\varepsilon,
\label{nkm}
\end{equation}
where $\varepsilon$ is a unit in $M$ and the variables are $X,Y\in\Z_M$.
This equation determines $X$ and $Y$ up to a unit factor in $M$.

Note that in our case, the fast algorithm \cite{relthue} for
finding "small" solutions of relative Thue equations is not applicable.
Also the reduction method corresponding to a possible application of
Baker's method (cf. \cite{book}) is also quite weak, 
producing relatively large reduced bounds,
therefore a huge amount of small values to test. 
These are the reasons which make the case
of sextic fields with a real quadratic subfield especially difficult.

In the case of real quadratic fields, $M$
has infinitely many units, and there are no tools to solve
 equation $J(\gamma)=1$  completely. Hence
we construct an efficient method to find "small" solutions 
(in the sense (\ref{bound}))
of the first factor, the relative Thue equation (\ref{nkm}), as well.

From equation (\ref{nkm}) it turns out that $X-\vartheta Y$ is a unit in $K$.
Denote by $\eta$ the fundamental unit of $M$, then in most cases
there are units $\varepsilon_1,\ldots,\varepsilon_h$ which, together with $\eta$
form a system of fundamental units of $K$. 
(For sextic fields with a real quadratic subfield we have $h=2,3,4$.)
By equation (\ref{nkm}) 
we may search for $X,Y$ in the form
\[
X-\vartheta Y=\zeta,
\]
with
\[
\zeta=\pm \eta^k \varepsilon_1^{k_1}\cdots \varepsilon_h^{k_h},
\]
where $k,k_1,\dots,k_h\in\Z$. 
If there are at least four conjugates $\zeta^{(i,j)}$, 
$(i,j)=(i_1,j_1)$, $(i_2,j_2)$, $(i_3,j_3)$, $(i_4,j_4)$,
with absolute value 
$<1$, then the system of equations
\begin{eqnarray}
x_1+\omega^{(i_1)} x_2 -\vartheta^{(i_1,j_1)} y_1-\omega^{(i_1)}\vartheta^{(i_1,j_1)} y_2
&=&\zeta^{(i_1,j_1)},\nonumber \\ 
\vdots &&    \label{sys3}\\
x_1+\omega^{(i_4)} x_2 -\vartheta^{(i_4,j_4)} y_1-\omega^{(i_4)}\vartheta^{(i_4,j_4)} y_2
&=&\zeta^{(i_4,j_4)} \nonumber
\end{eqnarray}
implies a very small bound for $x_1,x_2,y_1,y_2$, which values can be directly tested.
Hence we may assume that at least three conjugates 
$(i_1,j_1),(i_2,j_2),(i_3,j_3)$ of $\zeta$ are $>1$ in absolute value, that is 
$\log|\zeta^{(i_t,j_t)}|>0$, $t=1,2,3$.
For these conjugates we have
\[
|\log|\zeta^{(i_t,j_t)}||\le \log c_1,
\]
with 
\[
c_1=C+\ov{|\omega|}\ C+\ov{|\vartheta|}\  (C+\ov{|\omega|}\ C)\le c_0 C,
\]
where $c_0$ is a moderate constant and
$\ov{|\delta|}$ denotes the size of $\delta$ 
(the maximum absolute value of its conjugates).
$\zeta$ is a unit, hence 
the sum of all six $\log|\zeta^{(i,j)}|$ is zero.
Therefore the sum of the absolute values of the $\log|\zeta^{(i,j)}|$
with $\log|\zeta^{(i,j)}|<0$ is bounded by the sum of the absolute values
of the $\log|\zeta^{(i,j)}|$ with $\log|\zeta^{(i,j)}|>0$. Hence 
we have $|\log|\zeta^{(i,j)}||<3\log c_1$, if $\log|\zeta^{(i,j)}|<0$.

Consider
conjugates $(i_1,j_1)$, $\ldots$, $(i_{h+1},j_{h+1})$ such that 
in the system of equations in $k,k_1,\ldots,k_h$
\begin{equation}
\begin{array}{cccccc}
k\log|\eta^{(i_1)}|+&k_1\log|\varepsilon_1^{(i_1,j_1)}|+&\ldots& +k_h\log|\varepsilon_h^{(i_1,j_1)}|
&=&\log |\zeta^{^(i_1,j_1)}|, \\
\vdots &&&&&  \label{sys}\\
k\log|\eta^{(i_{h+1})}|+&k_1\log|\varepsilon_1^{(i_{h+1},j_{h+1})}|+&\ldots &+k_h\log|\varepsilon_h^{(i_{h+1},j_{h+1})}|
&=&\log |\zeta^{(i_{h+1},j_{h+1})}|, 
\end{array}
\end{equation}
the coefficient matrix is regular. 
(This is possible because the regulator of $K$ is non-zero.)
Solving this system of linear equations in  $k, k_1,\ldots,k_h$ by Cramer's rule 
we obtain an upper bound  for 
$\max(|k|,|k_1|,\ldots,|k_h|)\le B_0$ of magnitude $\log C$.
Note that we use this bound for $|k_1|,\ldots,|k_h|$, but 
 not for $k$, because 
we shall obtain a better bound for $k$ in Step 2.

Set $\xi=\eta^{-k}\zeta$, $X_0=x_{10}+\omega x_{20}=\eta^{-k}X$, 
$Y_0=y_{10}+\omega y_{20}=\eta^{-k}Y$.
The direct method would be to
enumerate all possible $k_1,\dots,k_h$ with absolute values $\le B_0$,
calculate 
$\xi=\pm\varepsilon_1^{k_1}\cdots \varepsilon_h^{k_h}$ and
solve the system of linear equations
\begin{eqnarray}
x_{10}+\omega^{(i_1)} x_{20} 
-\vartheta^{(i_1,j_1)} y_{10}-\omega^{(i_1)}\vartheta^{(i_1,j_1)} y_{20}
&=&\xi^{(i_1,j_1)},\nonumber\\
\vdots &&  \label{sys2}\\
x_{10}+\omega^{(i_4)} x_{20} 
-\vartheta^{(i_4,j_4)} y_{10}-\omega^{(i_4)}\vartheta^{(i_4,j_4)} y_{20}
&=&\xi^{(i_4,j_4)},\nonumber
\end{eqnarray}
in the variables $x_{10},x_{20},y_{10},y_{20}$
(here $(i_1,j_1),\ldots,(i_4,j_4)$ are distinct conjugates, 
so that the above system is uniquely solvable).
This would yield about $(2\log C)^h$ steps.

We can considerably diminish the number of $k_1,\dots,k_h$ to consider
by using a sieve method (see \cite{book}). 
Actually, this is one of the observations 
that make our method efficient.

Let $p$ be a prime such that 
the defining polynomial $f_{\Q}(x)=f^{(1)}(x)f^{(2)}(x)\in\Z[x]$ 
of $\alpha$ splits into linear
factors mod $p$:
\[
f_{\Q}(x)=(x-r_1)(x-r_2)(x-r_3)(x-r_4)(x-r_5)(x-r_6)\;\; (\bmod \; p).
\]
The conjugates $\alpha^{(i,j)}$ of the root $\alpha$ of $f_{\Q}(x)$ 
are then congruent
to one of $r_i$ modulo a prime ideal $\mathfrak p$, lying above $p$
in the ring of integers of a number field containing all 
conjugates. The $\theta,\varepsilon_1,\ldots,\varepsilon_h$ can
all be expressed in terms of $\alpha$, hence 
we can calculate 
integers $t_{(1,j)},e_{1,(1,j)}\ldots,e_{h,(1,j)}$ such that
\begin{equation}
\vartheta^{(1,j)}\equiv t_{(1,j)}\; (\bmod\; \mathfrak p), \;\;
\varepsilon_1^{(1,j)}\equiv e_{1,(1,j)}\; (\bmod\; \mathfrak p),\ldots,
\varepsilon_h^{(1,j)}\equiv e_{h,(1,j)}\; (\bmod\; \mathfrak p),
\label{ii}
\end{equation}
for $j=1,2,3$.
If we replace the algebraic numbers with these congruent integers 
in an identity, then the congruence modulo $\mathfrak p$ implies a 
congruence modulo $p$.

For this purpose apply Siegel's identity (cf. \cite{book}), say for $i=1$:
\[
(\vartheta^{(1,1)}-\vartheta^{(1,2)})(X-\vartheta^{(1,3)}Y)
\]
\begin{equation}
+
(\vartheta^{(1,2)}-\vartheta^{(1,3)})(X-\vartheta^{(1,1)}Y)+
(\vartheta^{(1,3)}-\vartheta^{(1,1)})(X-\vartheta^{(1,2)}Y)=0.
\label{siegel}
\end{equation}
We have $X-\vartheta^{(1,j)}Y=\eta^{(1)}\xi^{(1,j)}=
\pm (\eta^{(1)})^k
(\varepsilon_1^{(1,j)})^{k_1}\cdots (\varepsilon_h^{(1,j)})^{k_h}$.
We can simplify equation (\ref{siegel}) by $\eta^k$, then 
(\ref{ii}) implies
\[
(t_{(1,1)}-t_{(1,2)})e_{1,(1,3)}^{k_1}\cdots e_{h,(1,3)}^{k_h}
\]
\begin{equation}
+(t_{(1,2)}-t_{(1,3)})e_{1,(1,1)}^{k_1}\cdots e_{h,(1,1)}^{k_h}+
(t_{(1,3)}-t_{(1,1)})e_{1,(1,2)}^{k_1}\cdots e_{h,(1,2)}^{k_h}\equiv 0
\;\; (\bmod\; p).
\label{congr}
\end{equation}
We test this congruence for $|k_1|,\ldots,|k_h|\le B_0$.
Using integer arithmetics, 
this test is very fast, and their remains 
only a couple of exponents $k_1,\ldots,k_h$ that survive.

To emphasize the difference, note that 
in our Example 1 (see Section \ref{calc}), for $h=2$, $C=10^{50}$, 
to solve the system (\ref{sys2}) for all 
$k_1,\dots,k_h$ with absolute values $\le B_0$, took about 5 minutes,
using 250 digits accuracy of real numbers.
On the other hand, to test the congruence (\ref{congr}) for the same set of exponents
took only about 5 seconds. In that example 
out of the possibe 
93025 tuples $(k_1,\ldots,k_h)$  only 122 survived the congruence
test. 
For these  remaining exponent tuples we explicitly solved the system of equations
(\ref{sys2}) and  checked if the 
solutions $x_{10},x_{20},y_{10},y_{20}$ 
are integers. If so, we continued our calculation with 
$X_0=x_{10}+\omega x_{20}, Y_0=y_{10}+\omega y_{20}$. 
In Example 1 it took a negligible amount of time to solve
these remaining 122 systems of equations.
\\

\noindent
{\bf Step 2}\\
Next, we calculate generators of absolute power integral bases
(that is generators of power integral bases of $K$ over $\Q$).

Having calculated $X_0,Y_0$
we have (up to relative equivalence)
all generators $\gamma_0=\alpha X_0+\alpha^2 Y_0$ of relative power integral bases 
of $K$ over $M$.
Now we apply B of Lemma \ref{th2} and determine the unit factor $\nu$ of $M$
and $A\in\Z_M$, so that
$\gamma=A+\nu \gamma_0$
is a generator of a power integral basis of $K$.
Here  $A=a_1+\omega a_2\in\Z_M$, $\nu=\pm\eta^k$ is a unit in $M$.

We have 
\[
\gamma=a_1+a_2\omega\pm \eta^k\gamma_0=
a_1+a_2\omega\pm \eta^k (X_0\alpha+Y_0\alpha^2)=
a_1+a_2\omega\pm (\eta^k X_0)\alpha\pm (\eta^k Y_0) \alpha^2.
\]
According to (\ref{AXY}), we have
$X=x_1+\omega x_2=\pm \eta^k X_0,\; Y=y_1+\omega y_2=\pm\eta^k Y_0$.
Then 
\[
\gamma=a_1+a_2\omega+(x_1+\omega x_2)\alpha+(y_1+\omega y_2)\alpha^2.
\]
If $|k|$ is large, then $|x_1|,|x_2|,|y_1|,|y_2|$ also become large, therefore
(\ref{bound}) yields a bound also for $|k|$. 
The bound for $|k|$ that we can derive here is much better than the bound previously 
obtained for $|k|$ in Step 1.
More exactly,
if $X_0\ne 0$, then 
\begin{eqnarray}
x_1+\omega^{(1)} x_2&=&(\eta^{(1)})^k X_0^{(1)}, \nonumber \\
x_1+\omega^{(2)} x_2&=&(\eta^{(2)})^k X_0^{(2)}, \nonumber
\end{eqnarray}
which yields 
\[
x_2=\frac{(\eta^{(1)})^k X_0^{(1)}-(\eta^{(2)})^k X_0^{(2)}}{\omega^{(1)}-\omega^{(2)}}.
\]
If e.g. $|\eta^{(1)}|>1$ then $|\eta^{(2)}|<1$
and for $k>k_0$ we have  $|(\eta^{(2)})^k X_0^{(2)}|<0.1\cdot |\omega^{(1)}-\omega^{(2)}|$,
whence
\[
|x_2|>\frac{|\eta^{(1)}|^k |X_0^{(1)}|}{|\omega^{(1)}-\omega^{(2)}|}-0.1,
\]
and using (\ref{bound}) we obtain
\[
|\eta^{(1)}|^k< \frac{(C+0.1)\cdot |\omega^{(1)}-\omega^{(2)}|}{|X_0^{(1)}|},
\]
which gives a bound for $k$ from above of magnitude $\log C/\log |\eta^{(1)}|$.
A lower bound for $k$ can be obtained similarly. (The few small values of $k$ with absolute
value $\le k_0$ can be considered directly.)

To find the suitable $k$ and $a_2$ corresponding to $X_0,Y_0$, 
all we have to do is to run $k$ in the interval determined by the bounds obtained above,
for each $k$ to calculate $\eta^k (X_0\alpha+Y_0\alpha^2)$
and substitute $\gamma=a_2\omega\pm \eta^k (X_0\alpha+Y_0\alpha^2)$ into 
equation $J(\gamma)=1$ (cf. Lemma \ref{relindex}).
This equation is then a polynomial equation in $a_2$, whence it is possible
to find its integer roots (if any).

It is essential, that for a given $k$, the equation 
$P(a_2)=J(\gamma)-1=0$ is a polynomial equation of degree 9 with integer coefficients.
Even if for each $k$ 
we need multiple precision arithmetics to calculate the polynomial,
it is absolutely necessary to convert finally the coefficients to integers,
and (instead of calculating real roots of a polynomial with real coefficients),
calculate the integer roots of the polynomial with integer coefficients.
This is performed very fast.
This is a second point that makes our procedure efficient: in our Example 1
we had to consider values of the exponent $k$ with $|k|\le 120$.
For these values of $k$ we used 500 digits accuracy to calculate 
the polynomial $J(\gamma)$. To find the roots $a_2\in\Z$ of the equation 
$J(\gamma)=1$ for all these values of $k$,
using the polynomial with real coefficients took about 2-3 minutes.
On the other hand, converting the coefficients of the polynomial $J(\gamma)$
to integers and searching for integer roots $a_2$ of the equation $J(\gamma)=1$
took only 2-3 seconds.

\section{Examples and experiences}
\label{calc}

In this section first we give three examples for the application of our 
algorithm, then we formulate some general statements which were
inspired by the results of our numerical calculations.
\\

\noindent
{\bf Example 1.}\\
Let $M=\Q(\sqrt{2})$, and let $f(x)=x^3+2x+(1+\sqrt{2})$ be the relative 
defining polynomial of $\alpha$. The absolute defining polynomial
of $\alpha$ is $x^6+4x^4+2x^3+4x^2+4x-1$.
In $K=\Q(\alpha)$ a system of fundamental units is given by
\[
\eta=1+\sqrt{2},\;  \varepsilon_1=\alpha, \;
\varepsilon_2=-4+22\alpha-7\alpha^2+21\alpha^3-4\alpha^4+5\alpha^5.
\]
Define the bound $C=10^{50}$ for the coordinates of generators of
power integral bases with respect to the integral basis used above.\\
In (\ref{sys}) we obtained the upper bound $|k_1|,|k_2|\le 152$.
For the sieve method we took $p=809$, we had
\[
f_{\Q}(x)=x^6+4x^4+2x^3+4x^2+4x-1
\]
\[
\equiv (x+311)(x+36)(x+462)(x+536)(x+564)(x+518)\;\; (\bmod\; 809).
\] 
The congruence test of all $|k_1|,|k_2|\le 152$ took only a few seconds,
and only 122 pairs $k_1,k_2$ survived, out of $(2\cdot 152+1)^2=93025$ possible pairs.
For the remaining pairs $k_1,k_2$ we solved the system of equations 
(\ref{sys2}) and found three solutions:\\
$(k_1,k_2,x_{10},x_{20},y_{10},y_{20})=(0,0,1,0,0,0),(1,0,0,0,-1,0),(3,0,-1,-1,2,0)$.\\
To solve the system (\ref{sys2})  
for the remaining 122 pairs $k_1,k_2$
we used 250 digits accuracy and the calculation took a negligible amount of time.

Next, for each of the above solutions we calculated the upper bound for $|k|$,
as described in Step 2, which was usually about 120. 
Then for each $k$ we calculated 
\[
\gamma=a_2\omega\pm \eta^k (x_{10}+\omega x_{20}) \alpha
\pm \eta^k (y_{10}+\omega y_{20})\alpha^2,
\]
and substituted into $J(\gamma)=1$. 
To calculate the real coefficients of $J(\gamma)$ we used 500 digits accuracy.
Then we converted the coefficients of the polynomial $P(a_2)=J(\gamma)-1$
to integers and we were searching for the integer roots $a_2$ of $P(a_2)=0$.
This took only 2-3 seconds.

Solving this polynomial equation in $a_2$
we obtained the following integer solutions:\\
--for $(x_{10},x_{20},y_{10},y_{20})=(1,0,0,0)$, with $k=0$ we had $a_2=0$, resulting
$\gamma=\alpha$,\\
--for $(x_{10},x_{20},y_{10},y_{20})=(0,0,-1,0)$, with $k=1$ we had $a_2=-2$, resulting
\mbox{$\gamma=-2\sqrt{2}+\alpha^2(1-\sqrt{2})$},\\
--for $(x_{10},x_{20},y_{10},y_{20})=(-1,-1,2,0)$, we had no solutions.\\
Therefore this sextic field $K$ has, up to equivalence, the above two
generators of power integral bases, having coefficients $<10^{50}$ in the integral basis.

\vspace{0.5cm}

In the following two examples we performed similar calculations with $C=10^{50}$.
The fields have the same signature, the bounds we obtained and the
CPU times were also similar. 
We only list the solutions using the notation of Example 1.
\\

\noindent
{\bf Example 2.}\\
$f(x)=x^3+(2+\sqrt{2})x+1$\\
$(x_{10},x_{20},y_{10},y_{20})=(1,0,0,0)$, $k=0, a_2=0$,  $\gamma=\alpha$,\\
$(x_{10},x_{20},y_{10},y_{20})=(0,0,-1,0)$, $k=0, a_2=-1$,  $\gamma=-\sqrt{2}-\alpha^2$,\\
$(x_{10},x_{20},y_{10},y_{20})=(-1,0,2,0)$, no solutions for $a_2$. \\
\\

\noindent
{\bf Example 3.}\\
$f(x)=x^3+(2+\sqrt{2})x+(1+\sqrt{2})$\\
$(x_{10},x_{20},y_{10},y_{20})=(1,0,0,0)$, $k=0, a_2=0$,  $\gamma=\alpha$,\\
$(x_{10},x_{20},y_{10},y_{20})=(-2,2,-3,1)$, no solutions for $a_2$,\\
$(x_{10},x_{20},y_{10},y_{20})=(0,0,-1,1)$, $k=0, a_2=0$,  $\gamma=\alpha^2(-1+\sqrt{2})$,\\
$(x_{10},x_{20},y_{10},y_{20})=(3,-2,-6,4)$, $k=2, a_2=-4$,  $\gamma=-4\sqrt{2}+\alpha-2\alpha^2$.\\

We performed calculations for several sextic fields with a real quadratic subfield.
Among others we considered relative defining polynomials of type
$f(x)=x^3+(2+b\sqrt{2})x+(1+c\sqrt{2})$ of $\alpha$
over $M=\Q(\sqrt{2})$, with some values of $b$ and with $c=0,\pm 1$.
These $\alpha$ have absolute defining polynomial 
\[
g(x)=x^6 + 0\cdot x^5+ 4x^4 + 2x^3 + (4-2b^2)x^2 + (-4bc + 4)x - 2c^2 + 1.
\]
We observed, that in the number field generated by the root $\alpha$ of these monogenic
polynomials, in addition to $\alpha$, there is another generator of
power integral basis, of the form
\[
\gamma=(4-2b^2)\alpha+2\alpha^2+4\alpha^3+0\cdot \alpha^4+\alpha^5.
\]
The coordinates of $\gamma$ are just the coefficients of the polynomial $g(x)$
in the reverse order. Further, for $c=0,\pm 1$ the constant term of $g(x)$
is $\pm 1$, hence $\alpha$ is a unit. 
This leads us to the following statements, which is easy to prove, 
but can be useful in some cases.

\vspace{1cm}

\noindent
{\bf Statement.} {\it If a unit $\alpha$ generates a power integral basis 
in the number field $K=\Q(\alpha)$, then $1/\alpha$ also generates a power integral
basis in $K$.}

\vspace{1cm}

\noindent
{\bf Proof.} Denote by $\alpha^{(i)},(i=1,\ldots,n)$ the conjugates of $\alpha$.
Then
\[
D\left(\frac{1}{\alpha}\right)=\prod_{1\le i<j\le n}
\left(\frac{1}{\alpha^{(i)}}-\frac{1}{\alpha^{(j)}}\right)^2
=
\prod_{1\le i<j\le n}\frac{(\alpha^{(i)}-\alpha^{(j)})^2}{(\alpha^{(i)}\alpha^{(j)})^2}
=\frac{D(\alpha)}{N(\alpha)^{2(n-1)}}=D(\alpha),
\]
since $N(\alpha)=\pm 1$.\mbox{}\hfill $\Box$

\vspace{1cm}

\noindent
{\bf Corollary.} {\it Assume the number field $K$ is generated by a root 
$\alpha$ of an irreducible monogenic polynomial 
\[
f(x)=x^n+a_{n-1}x^{n-1}+\ldots +a_2x^2+a_1x\pm 1\in\Z[x],
\]
of degree $>2$ with constant term $\pm 1$.
Then, in addition to $\alpha$, the inequivalent element
\[
a_2\alpha+a_3\alpha^2+ \ldots + a_{n-1}\alpha^{n-2}+\alpha^{n-1}
\]
also generates a power integral basis in $K$}.

\vspace{1cm}
\noindent
{\bf Proof.} 
Monogenity of the polynomial $f(x)$ means, that $\alpha$ generates a power integral
basis in $K$. The constant term of the polynomial is $\pm 1$, hence
$\alpha$ is a unit. Then, by above statement $1/\alpha$ also generates 
a power integral basis in $K$. 
By
$\alpha^n+a_{n-1}\alpha^{n-1}+\ldots +a_2\alpha^2+a_1\alpha\pm 1=0$, we obtain
\[
\pm \frac{1}{\alpha}=-a_1-a_2\alpha-\ldots -a_{n-1}\alpha^{n-2}-\alpha^{n-1},
\]
which is equivalent to the element in the Corollary..

If  $1/\alpha$  were equivalent to $\alpha$, then $\alpha$ would be a root
of a second degree polynomial. Therefore if the degree of $f(x)$ is
$>2$, then $\alpha$ and $1/\alpha$ are inequivalent elements.\mbox{}\hfill $\Box$

\vspace{1cm}

\noindent
{\bf Acknowledgements} The author is thankful to the referee for the 
careful reading and for several useful remarks that contributed
to improve the quality fo the paper.

\end{document}